\documentstyle[amssymb]{amsart}

\numberwithin{equation}{section}

\newtheorem{theorem}{Theorem}[section]
\newtheorem{lemma}[theorem]{Lemma}
\newtheorem{proposition}[theorem]{Proposition}

\theoremstyle{definition}

\def\R{{\Bbb R}}

\def\Z{{\Bbb Z}}
\def\C{{\Bbb C}}
\def\P{{I \! \! P}}
\def\D{\hbox{\goth D}}
\def\End{\operatorname{End}}

\def\nn#1{}
\def\Ca{\C[ x^\alpha]}

\begin{document}
\title[Differential operators and invariant monomials]{Classification
of linear differential operators with an invariant subspace of monomials}

\author{Gerhard Post}
\address{Gerhard Post, Department of Applied Mathematics\\
 University of Twente\\ P.O. Box 217\\
 7500 AE Enschede, The Netherlands}
\email{post@@math.utwente.nl}
\author{Alexander Turbiner}
\thanks{A.T. is supported in part
by CAST grant, US National Academy of Sciences}
\address{Alexander Turbiner, Institute for Theoretic and Experimental
Physics, \newline 117259 Moscow, Russia (on leave in absence)}
\curraddr{Alexander Turbiner, Mathematics Department,
Case Western Reserve University, Cleveland OH 44106}
\email{turbiner@@vxcern.cern.ch}
\date{July 9, 1993}
\maketitle

\begin{abstract}
A complete classification of linear
differential operators possessing finite-dimensional invariant subspace
with a basis of monomials is presented.
\end{abstract}

\section{Introduction}

One of the old-standing problems in the theory of special functions is the
classification of all linear differential operators that admit an infinite
sequence of orthogonal eigenvectors in the form of polynomials.
See \cite{L} for an overview.
In 1929 S. Bochner \cite{B} had solved this problem for second order
differential operators on the complex (real) line. Therefore we
name this problem the ``Bochner problem''.

The main purpose of the present paper is to study a more general
problem, which we name (following \cite{T2}) the ``generalized Bochner
problem''. We ask for a classification of all linear differential
operators, which possess a {\bf finite-dimensional} invariant subspace
of polynomials. It turns out this problem is rather sophisticated.
However, things get much more tangible when, instead of invariant
subspace in polynomials, we require that the (finite-dimensional)
invariant subspace has a basis of monomials. This problem is solved
completely (Section 3). In \cite{T2} a particular case of this problem
was solved, namely if the invariant subspace is the linear space of
polynomials of degree not higher than some fixed integer.
 The classification of the linear
operators possessing such an invariant subspace was given through the
universal enveloping algebra of the algebra $sl_2$ taken in a special
representation by first-order differential operators.

The results of Section 3 have an impact to the general problem. This
is presented in Section 4.  In Section 5, we give the explicit expressions for
the
second-order differential operators $T_2$, possessing a
finite-dimensional invariant subspace in monomials. These
operators are of great interest for finding explicit solutions to the
Schr\"odinger equation
\begin{equation}\label{eq1}
(-\frac{d^2}{dx^2} + V(x))\Psi(x) = \lambda \Psi(x)
\end{equation}
since the eigenvalue problem for the operator $T_2, T_2\varphi=\lambda
\varphi$ can be reduced to the Schr\"odinger equation by a change of the
variable, $x'=x'(x)$, and introducing a new (gauge-transformed)
function, $\Psi=\varphi\exp(-a(x))$.
Hence the operators described in Section 5 lead to a special class of
quasi-exactly-solvable Schr\"odinger equations complementary those
described already \cite{T2}.

Finally we make some concluding remarks in Section 6 concerning the
case where the powers of $x$ are not natural numbers.

Though in this paper we work over the complex numbers, the main
results hold for real numbers as well.

\section{Generalities}

Here we state some general results concerning differential operators
which leave  a finite-dimensional subspace of $\C [x]$
invariant. The results of this section can easily be extended to more
variables.
Let \D \ denote the algebra of linear (finite-order)
differential operators on $\C [x]$ with polynomial coefficients. We
denote the symbol $\frac{\displaystyle d}{\displaystyle dx}$ by $\partial$.

Let $V$ be a finite-dimensional subspace in $\C [x]$ of dimension $n$,
and let $T: \C [x] \to \C [x]$ be a linear operator. It is easy to
show that any linear operator $T$ can be written as an infinite
order linear differential operator:
\[
	T = \sum_{i =0}^\infty P_i \partial^i,
\qquad P_i \in \C [x].
\]

This can be proved easily by an inductive construction:
\[
	P_0 = T(1), \quad P_1 = T(x) - xP_0, \ldots .
\]
It follows immediately that the action of $T$ on $V$ can be
represented by a finite-order differential operator, since $V \subset
\P_k$ (where $\P_k$ denotes the space of polynomials of degree not higher than
$k$, and $k$ sufficiently large), and hence for $T_k$:
\[
	T_k = \sum_{i=0}^k P_i \partial^i ,
\]
we have that $T_k(v) = T(v)$ for all $v \in V$.

So the following proposition is immediate:

\begin{proposition}
Let $V \subset \C [x]$ be a subspace of dimension $n$, and let $\D_V$
denote the algebra of differential operators that
leave $V$ invariant (i.e. $T_k \in \D_V \Rightarrow T_k(V) \subset V$).
Then $\D_V$ is isomorphic to the semi-direct product of $\End (V)$ and
$I$, where $\End (V)$ is the algebra of linear operators on $V$, and
the ideal $I$ is the algebra of differential operators which
annihilate $V$. \hfill{$\Box$}
\end{proposition}

We can paraphrase this proposition in the following way. As
representation on $V$, $\D_V$ is the full $n \times n$-matrix algebra, and the
kernel of this representation is an infinite dimensional ideal in
$\D_V$.

Another consequence of this observation is, that for
any element $v \in \C [x]$
and $\lambda \in {\bold C}$, we can find a differential operator
$T$ such that $T v=\lambda v$.

In this paper, we mainly discuss the case that $V$ is
graded. Now $\C [x]$ is a graded algebra by putting $\deg (x^k) = k$.
Let us take $V$ a graded subspace. This means exactly that $V$ has a
basis of monomials. Hence we assume that $V = \langle x^{i_1},x^{i_2},
\ldots , x^{i_n} \rangle$, which we abbreviate to $V = \langle
x^I \rangle$, $I = \{i_1,i_2,\ldots , i_n\}$.

The fact that $\C [x]$ is graded has as a consequence that $\End(\C [x])$
is also graded, putting $\deg(T) = m$, if $T(x^i) \in \langle x^{i+m}
\rangle$ for all $i$. Moreover for graded $V$ it follows that $\D_V$ is
also graded. So, in order to describe the structure of $\D_V$, it is sufficient
to
describe the homogeneous components of $\D_V$. This is the concern
of the next section.

\section{Algebras leaving a space of monomials invariant}

As before, we are interested in finite-order linear differential
operators $T$, such that $T(V) \subset V$. We assume $T$ to be graded
of degree $m$ and order $k$, which means that
\[
	T = \sum^k_{i=0} c_ix^{i+m} \partial^i.
\]
Here $c_k \not = 0$ and $m \in \Z$. Moreover $c_i = 0$ for $i+m < 0$.
In particular, we see that if the degree of $T$ is negative, say
deg$(T)=-m, m>0$, then the order of $T$ is at least $m$.

The following lemma plays a crucial role in our classification.

\begin{lemma}
Let $T$ be a differential operator of degree $m$ and order $k$.
\begin{enumerate}
\item Suppose $m \geq 0$. Then there exist numbers
$\alpha_1,\ldots,\alpha_k \in \C$ such that
\[
	T = c_kx^m (x\partial - \alpha_1)(x\partial -\alpha_2)\ldots
	(x \partial - \alpha_{k}).
\]
\item Suppose $m < 0$, hence $k \geq -m$. Then there exist numbers
$\alpha_1,\ldots,\alpha_{k+m} \in \C$ such that
\[
	T = c_k\partial^{-m} (x\partial - \alpha_1)(x\partial -\alpha_2)\ldots
	(x \partial - \alpha_{k+m}).
\]
\end{enumerate}
\end{lemma}
\begin{pf}
If $m \geq 0$ it is clear that any differential operator of degree $m$ is of
the
form $T = x^m \sum\limits_{i=0}^k c_ix^i\partial^i$. Similarly, if
$m<0$, $T$ can be put in the form $T = \partial^{-m}
\sum\limits_{i=0}^k \tilde{c}_ix^i\partial^i$ with $\tilde{c}_k=c_k$. So it
suffices to prove the
lemma for $m=0$. But then we have $T(x^\alpha) = P_k(\alpha)x^\alpha$,
where $P_k$ is a $k$-th order polynomial. So $P_k(\alpha) =
c_k(\alpha - \alpha_1)(\alpha -
\alpha_2) \cdots (\alpha -\alpha_k)$ for some $\alpha_i \in \C$. Now
\[
	((x\partial - \alpha_1) \cdots (x\partial -
	\alpha_k))(x^\alpha) = (\alpha - \alpha_1)(\alpha - \alpha_2)\cdots
	(\alpha - \alpha_k),
\]
{}From this it follows that $T = c_k(x\partial - \alpha_1)(x\partial -
\alpha_2) \cdots (x\partial - \alpha_k)$, since the representation of \D \
on $\C[x]$ is faithful.
\end{pf}
It is easy to see that this factorization is unique up to the ordering of the
factors
$x\partial - \alpha_i$. However, a certain ordering of the
factors has no meaning, since these factors commute.

Above representation of $T$ given by lemma 3.1 is very convenient, when
we study the operators $T$ which leave $V = \langle x^{i_1},\ldots,
x^{i_n}\rangle$ invariant, i.e. $T \in \D_V$. We introduce the
following notation. For $I = \{i_1,i_2,\ldots,i_n\}$ and $m \in \Z$,
we put
\[
	I^{(m)} = \{ i \in I \mid i+m \geq 0 {\text{ \ and }} i+m \not\in I\}.
\]
If $\deg (T) = m$ and $T \in \D_V$, then it is clear that with necessity
$T(x^i) = 0$ for $i \in I^{(m)}$, since $T(x^i) = c\cdot x^{i+m}$, but
$x^{i+m} \not\in V$, so $c = 0$. This leads to the following
\begin{theorem}
Let $V = \langle x^{i_1},\ldots, x^{i_n} \rangle$ and $T$ a
finite-order differential operator such that $\deg (T) = m$. Suppose
$I^{(m)} = \{\alpha_1,\ldots,\alpha_k\}$. Then $T \in \D_V$, if and only
if
\[
	T = \widetilde{T} \cdot (x\partial - \alpha_1)(x\partial -
	\alpha_2) \cdots (x\partial - \alpha_k),
\]
where $\widetilde{T}$ is some differential operator of degree $m$.
\end{theorem}
\begin{pf}
The if-part is trivial. So assume that $T \in \D_V$, and suppose it has
order $s$. According to lemma 3.1, for $m \geq 0$, $T$ can be represented in
the form
\[
	T = c \cdot x^m (x\partial - \beta_1)(x\partial -
	\beta_2) \cdots (x\partial - \beta_s) \qquad (c \in \C, c \neq 0).
\]
We need that $T(x^i) = 0$ for $i \in I^{(m)}$. On the other hand,
we have
\[
	T(x^i) = c(i-\beta_1)(i-\beta_2)\cdots(i-\beta_s)x^{i+m}.
\]
Hence it follows that $\{\alpha_1,\ldots,\alpha_k\} \subset
\{\beta_1,\ldots,\beta_s\}$. After rearranging the $\beta$'s we find
\[
	T = c\cdot x^m(x\partial -\beta_1)(x\partial -\beta_2) \cdots
(x\partial - \beta_{s-k})(x\partial -\alpha_1)\cdots(x\partial -
\alpha_k).
\]
So for $m \geq 0$ the proposition is proved; $\widetilde{T} =c\cdot
x^m(x\partial -\beta_1)(x\partial -\beta_2) \cdots (x\partial -
\beta_{s-k})$. For $m<0$ the proof is similar.
\end{pf}

{\bf REMARK.} From the previous
proposition it follows that the order of $T \in \D_V$ with $\deg(T) =
m$ is at least $k$, where $k$ is the number of elements in $I^{(m)}$.
In fact, up to a scalar coefficient, the element of order $k$ is
unique:
\[
	T = x^m(x\partial - \alpha_1) \cdots (x\partial - \alpha_k),
\]
where $\{\alpha_1,\ldots,\alpha_k\} = I^{(m)}$.

If $m=0$, $I^{(m)} = \emptyset$, and we have $T=1$. So
all differential operators of degree 0 are in $\D_V$; these elements form a
commutative subalgebra of $\D_V$ generated by $x\partial$.

\section{The case of polynomial subspaces}
In section 3, we performed the classification of differential
operators with an invariant subspace $V$ that has a basis of monomials.
If $V$ has no basis of monomials, $\D_V$ is not be graded,
but only filtered. This causes a major difficulty. However, considering the
corresponding grading, we still can deduce some properties of $\D_V$ in
this case.

So, let $V$ have a basis of the form
\[ x^{i_1}+c_{11}x^{i_1-1}+\dots ,\quad x^{i_2}+c_{21}x^{i_2-1}+\dots ,\quad
\dots\quad  ,\quad x^{i_n}+c_{n1}x^{i_n-1}+\dots \]
We can assume that all $i_j$ are different. The graded space $V(g)$
associated to $V$ is $\langle x^{i_1},\ldots, x^{i_n}\rangle$.

Let $T\in \D$ of order $k$ be of the form
\[ T = \sum_{i=-k}^{m} T^{(i)} \]
with $T^{(i)}$ of degree $i$, and $T^{(m)} \not = 0$. (An operator of
order $k$ has degree $-k$ or higher, the term of degree $-k$ being a
multiple of $\partial^k$). We call $T^{(m)}$ the associated
graded operator. Now it is easy to prove
\begin{theorem}
Let $V,\ V(g),\ T$ and $T^{(m)}$ be as above.\\
 If $T \in \D_V$, then $T^{(m)} \in \D_{V(g)}$.
\end{theorem}
\begin{pf}
Suppose $T \in \D_V$, and consider $T(x^{i_j}+c_{j1}x^{i_j-1}+\dots)$.
If $m=0$ there is nothing to prove, so consider $m \neq 0$. Let $I =
\{i_1,i_2,\ldots,i_n\}$. If $i_j \in I^{(m)}$ then $T^{(m)}(x^{i_j})$
should be 0, since no term of $T$ can cancel this term. Hence we find
exactly that $T^{(m)}(x^{i_j})=0$ for all $i_j \in I^{(m)}$, i.e.
$T^{(m)} \in \D_{V(g)}$.
\end{pf}

{}From this we derive an easy corollary:\\
{\it If $T_k \in \D$ of order $k$ possesses an infinite number of
(linearly independent)
eigenvectors, then $T = \sum_{i=-k}^{0} T^{(i)}$, degree of $T^{(i)}= i$
and $T^{(0)} \neq 0$. }

A similar reasoning can be performed for the part of $T$ with minimal
degree, but this seems to give not much information.

\section{Classification of second order differential operators}
\subsection{Generic situation}
Now we proceed to 2-nd order differential operators
$T_2$, which admit a finite-dimensional invariant space of polynomials
with a basis of monomials.  We are interested in this problem in
connection with finding explicit solutions to the Schr\"odinger
equation (1.1). As mentioned in the Introduction, this involves some
transformations see \cite{T1}. These transformations are of 2 types:
the first is a change of basis, and the
second is called  ``gauge'' transformation, which amounts to
changing $T$ to $g T g^{-1}$, where $g$ is a non-zero function.
To delete some ambiguity in our spaces of
monomials, we impose 2 conditions on $V = \{x^{i_1},x^{i_2},\ldots,x^{i_n}\}$:
\begin{enumerate}
\item We assume that $1 \in V$ (it removes an ambiguity resulting from
gauge transformations).
\item We assume that gcd$(i_1,i_2,\ldots,i_n)= 1$, i.e. that the
powers have no common factor (it removes an ambiguity resulting from
changes of variable).
\end{enumerate}
For the classification of differential operators in $\D_V$ these
assumptions do not make much difference. If $T_2= \sum T^{(i)}$, where the
degree
of $T^{(i)}$ is i, these two assumptions effect only the terms
$T^{(-1)}$ and  $T^{(-2)}$.

So let us start the classification. Suppose $T_2 \in \D_V$, $T_2=\sum_m
T^{(m)}$
where $T^{(m)}$ has degree $m$ and order less or equal to 2. If
$T^{(m)}$ is non-zero, then according to proposition 3.2,  $I^{(m)}$
contains 0,1 or 2 elements.
We distinguish these 4 cases:
\begin{enumerate}
\item
For all $m>0$, $I^{(m)}$ contains more than 2 elements. So $T_2$
contains no terms of positive degree. Hence
$T_2$ preserves $\langle 1,x,x^2,\ldots,x^{n}\rangle$ for all $n$.
This type of operators is already studied in \cite{T2} called there
exactly-solvable operators), and we do not repeat it here.
\item
$I^{(m)}$ is empty. This can only happen, if $m=0$. But this case
is trivial, since all operators of degree 0 are in $\D_V$. Hence
degree 0 contributes to $T_2$ the operator
$\alpha_1x^2\partial^2+\alpha_2x\partial+\alpha_3$. This part we call
trivial, and is always present in $T_2$.
\item
$I^{(m)}$ contains one element.
Let us assume here that $i_1>i_2>\dots >i_n$. Then we have
$i_j=i_{j-1}+m$ so that $i_j =(n-j)m$ (since we assumed that $1 \in
V$, so $i_n=0$). But we also assumed that the $i_j$ have no common
factor, so it follows that $m=1$, and hence $V=\langle
1,x,x^2,\ldots,x^{n-1}\rangle$. Hence we are in the case that is
extensively discussed in \cite{T2}.
\item
$I^{(m)}$ contains 2 elements, $m>0$, and no $I^{(l)}$
for $l>0$ contains one element. Suppose $i_1$ and $i_2$ are
the 2 elements in $I^{(m)}$.
Then the set $\{ i_1,i_2,\ldots,i_n \}$ is of the following
specific form:
\[ i_3=i_1-m,\ i_5=i_3-m,\ \ldots ,\ i_{2r-1}=i_{2r-3}-m \]
and
\[ i_4=i_2-m,\ i_6=i_4-m,\ \ldots ,\ i_{2s}=i_{2s-2}-m \]
We call $i_1,i_3,\ldots,i_{2r-1}$ and $i_2,i_4,\ldots ,i_{2s}$ {\it
chains} with step $m$ and length $r$ and $s$, respectively. In
general, with no special relation, by which we mean that $I$ can be split
into two chains in exactly one way, the most general second-order
operator in $\D_V$ is
\[ T_2 = \alpha_1x^m(x\partial-i_1)(x\partial-i_2) + \alpha_2x^2\partial^2
+ \alpha_3x\partial + \alpha_4 \]
and $T_2$ contains an extra term in two cases:
\begin{itemize}
\item[(a)] If $m=1$, (and therefore $i_{2s}=0$), we get an extra term
\linebreak $\alpha_5\partial(x\partial-i_{2r-1})$.
\item[(b)] If m=2 and $\{i_{2s},i_{2r-1}\} = \{0,1\}$, we get an extra
term $\alpha_5\partial^2$.
\end{itemize}
All this can easily be proved, by examining the possibilities for
which $I^{(m)}$ can have 0,1 or 2 elements.
\end{enumerate}

\subsection{Special subspaces}
For existence of non-trivial second order operators in $\D_V$, it is
necessary that $I$ can be split into 2 chains. The form of $T_2$ above
is under the assumption that $I$ can be split into two chains in
exactly one way. There are cases, where $I$ can be split in more
than one way. As an example consider $I=\{0,1,2,\ldots,98,100\}$. Then
$I=\{0,1,2,\ldots,98\} \cup \{100\}$ or $I=\{0,2,4,\ldots,98,100\} \cup
\{1,3,5,\ldots,97\}$. Consequently the space $V$ admits a more general
second-order operator. Here we describe all special cases, which
fall into 4 groups:
\begin{itemize}
\item[{\bf Case A.}]
The dimension of $V$ is 3, so $I=\{0,m,m+l\}$, $l \neq m$. Then
\begin{align*}
T_2 & = \alpha_1 x^{l+m}(x\partial -m)(x\partial -l-m) + \alpha_2
x^{l+1}\partial(x\partial -l-m)\\
& \quad + \alpha_3 x^m(x\partial -m)(x\partial -l-m) + \alpha_4
x^2\partial^2 + \alpha_5 x\partial +\alpha_6
\end{align*}
$T_2$ gets a certain extra terms in the following cases:
\begin{itemize}
\item[(a)]
$m=1$, $l > 2$. Extra term $\alpha_7\partial(x\partial -l-m)$.
\item[(b)]
$m=1$, $l = 2$. Extra terms $\alpha_7 \partial(x\partial
-3) + \alpha_8 \partial^2$.
\item[(c)]
$l=1$ (hence $m>1$). Extra term $\alpha_7\partial(x\partial -m)$.
\end{itemize}
\item[{\bf Case B.}]
The dimension of $V$ is 4, and $I$ is ``symmetric'', i.e.
$I=\{0,m,m+l,2m+l\}$, $l \neq m$. Then
\begin{align*}
T_2 & = \alpha_1 x^{l+m}(x\partial -2m-l)(x\partial -l-m) + \alpha_2
x^m(x\partial -2m-l)(x\partial -m)\\
& \quad + \alpha_3 x^2\partial^2 + \alpha_4
x\partial + \alpha_5
\end{align*}
$T_2$ gets an extra term only if
\begin{itemize}
\item[(a)]
$m=1$. Extra term $\alpha_6\partial(x\partial -m-l)$.
\end{itemize}
\item[{\bf Case C.}]
$I$ has one runner ahead at distance $2$, i.e. $I=\{0,1,2,\ldots,n-2,n\}$.
Then $T_2$ takes the form
\begin{align*}
T_2 & = \alpha_1 x^{2}(x\partial -n)(x\partial -n+3) + \alpha_2
x(x\partial -n)(x\partial -n+2)\\
& \quad + \alpha_3 x^2\partial^2 + \alpha_4x\partial + \alpha_5\\
& \quad + \alpha_6\partial(x\partial-n) + \alpha_7\partial^2
\end{align*}
\item[{\bf Case D.}]
$I$ has one runner left behind at distance $2$, i.e. $I=\{0,2,3,4,\ldots,n\}$.
Then $T_2$ is of the form
\begin{align*}
T_2 & = \alpha_1 x^{2}(x\partial -n)(x\partial -n+1) + \alpha_2
x^2\partial(x\partial -n)\\
& \quad + \alpha_3 x^2\partial^2 + \alpha_4
x\partial + \alpha_5 \\
& \quad + \alpha_6\partial(x\partial - 2)
\end{align*}
\end{itemize}

These special cases exhaust the list of all exceptional cases, which do
not belong to the general classification given above. Let us
prove this. It is clear that we can assume that the dimension is 5 or
more, since the dimensions 3 and 4 are considered above. We know that
$I$ can be split into two chains of step, say, $m$, so $ i_3=i_1-m,\
i_5=i_3-m,\ \ldots ,\ i_{2r-1}$ and $i_4=i_2-m,\ i_6=i_4-m,\ \ldots ,\
i_{2s}$. Suppose that $I$ can also be split in two chains of step $l$.
We may assume that $l>m$. Now we distinguish three cases:
\begin{enumerate}
\item
$r=1$, so $i_3$ is not present.\\
We have either $I^{(l)} \supset \{i_1,i_2\}$ or $I^{(l)} \supset
\{i_1,i_4\}$. Therefore we consider two subcases:
\begin{itemize}
\item[a.] $I^{(l)} \supset \{i_1,i_2\}$, so $l \neq i_1-i_2$. Since
$\mid I^{(l)} \mid = 2$, we have $i_4 \not \in I^{(l)}$, so $i_1-i_4
= l$. But then $i_1-i_6 \neq l$, and therefore $i_2-i_6 = 2m = l$.
This is not allowed, since $i_1-i_4 = 2m$ implies that $i_1-i_2 = m$.
\item[b.] $I^{(l)} \supset \{i_1,i_4\}$, so $i_1-i_2 = l$. But also
$i_6 \not \in I^{(l)}$, and this is only possible if $i_2-i_6 = l=2m$.
Therefore this configuration leads to case C.
\end{itemize}
\item
$s=1$, so $i_4$ is not present.\\
Clearly $I^{(l)} \supset \{i_1,i_3\}$. Moreover $i_2 < i_3$ since
otherwise also $i_2 \in I^{(l)}$. But then $i_5 \not \in I^{(l)}$
implies that $i_5+l=i_1$, and hence $l=2m$. Again using $i_2 \not \in
I^{(l)}$ implies that $i_2+l=i_{2r-1}$. Therefore this configuration
leads to case D.
\item
$r >1$ and $s>1$, so $i_3$ and $i_4$ are both present.\\
Always $I^{(l)} \supset \{i_1,i_3\}$, and necessarily
$i_2 < i_3$, since otherwise also $i_2 \in I^{(l)}$.  Again we have
two subcases:
\begin{itemize}
\item[a.]
$i_5$ is present. We need $i_5+l=i_1$, so $l=2m$. Like in case 2
above, it follows that $i_2= i_{2r-1}-2m$. But then $i_4 \in I^{(l)}$.
\item[b.]
$i_6$ is present, but not $i_5$. We need either $i_1-i_2 = l$ or
$i_3-i_2 = l$.

If $i_1-i_2=l$, then $i_6+l=i_2$ is the only possibility, so $l=2m$,
and hence $i_3-i_2=m$. This is forbidden because then we have one chain.

If $i_3-i_2 =l$, then $i_4+l \not\in \{i_1,i_3,i_2\}$, so $i_4 \in
I^{(l)}$. So this is also impossible.
\end{itemize}
\end{enumerate}

\section{Conclusion}

In the previous sections we considered the case that the powers of
$x$ are natural numbers. From algebraic point of view, this is not
a crucial restriction. One could take the powers to be all
integer, rational, real or complex numbers, or any other abelian
subgroup of $\C$. (If one takes $\C$, one has to define a suitable
total ordering on $\C$ to be able to consider positive and negative).
This has as a
main advantage that such algebras of ``polynomials'' admit the
isomorphism $x \mapsto 1/x$. Moreover, the map
$x^k \mapsto x^{k+l}$ is a linear isomorphism; it is the gauge
transformation discussed before.

An even more special case is that the set of allowed powers form a
field. In this case, the change of basis $x\mapsto x^{m}$ is an
isomorphism.

We shortly discuss the case that the powers of $x$ are real numbers.
The algebra of ``generalized polynomials'', which we denote
(suggestively) by $\Ca$ is as a linear space
\[ \bigoplus_{\alpha \in \R} \C \cdot x^\alpha \]
and the multiplication is given by
\[
x^\beta \cdot x^\gamma = x^{\beta+\gamma}\qquad \beta, \gamma \in \R
\]
Let $\D$ denote the algebra of differential operators with
coefficients in $\Ca$ now. The degree of $T \in \D$ can be any
real number. For $T \in \D$, deg$(T)=m$ we have
\[
	T = \sum^k_{i=0} c_ix^{i+m} \partial^i.
\]
with $c_i \in \C$, and $c_k \neq 0$. Lemma 3.1 for this operator $T$
now looks like:
\begin{lemma}
 There exist numbers
$\alpha_1,\ldots,\alpha_k \in \C$ such that
\[
	T = c_kx^m (x\partial - \alpha_1)(x\partial -\alpha_2)\ldots
	(x \partial - \alpha_{k}).
\]
\end{lemma}
In particular, the case that $m<0$ dissappears, or, better,
$\partial^{-m}$ factorizes:
\[
\partial^{s} = x^{-s}(x \partial - (s-1))(x \partial - (s-2)) \cdots (x
\partial - 1)x \partial \qquad (s > 0)
\]

We denote $I=\{ i_1,i_2,\ldots,i_n\}$, where $i_j \in \R,\
j=1,\ldots,n$, and define for $m \in \R$
\[
	I^{(m)} = \{i \in I \mid i+m \not\in I\}.
\]
Let $V=\langle x^{i_1},x^{i_2},\ldots,x^{i_n}\rangle$ and let $\D_V$ denote
the algebra of operators in $\D$ that leave $V$ invariant. One can
easily check that Theorem 3.2 holds literally.
Due to the remarks at the beginning of this section, we have the following.
\begin{proposition}\hfill {}
\begin{enumerate}
\item The gauge transformation  $T \mapsto x^lTx^{-l}$ gives an
isomorphism between $\D_V$ and $\D_W$, with $V=\langle
x^{i_1},x^{i_1},\ldots,x^{i_n}\rangle$ and $W=\langle
x^{i_1-l},x^{i_2-l},\ldots,x^{i_n-l}\rangle$.
\item The change of basis $x'=x^m$ induces an isomorphism between
$\D_V$ and $\D_W$, where $V=\langle x^{i_1},x^{i_1},\ldots,x^{i_n}\rangle$ and
$W=\langle x^{i_1/m},x^{i_2/m},\ldots,x^{i_n/m}\rangle$.
\end{enumerate}
\end{proposition}
\newpage
{\it Proof.}
\begin{enumerate}
\item Part (1) is obvious.
\item By the change of basis $x'=x^m$ we have that
\[ x\partial = m\; x'\partial', \qquad \partial'=\frac{d}{dx'}.  \]
Hence the operator $T = \tilde{T}\cdot (x\partial - \alpha_1)(x\partial -
	\alpha_2) \cdots (x\partial - \alpha_k)$, where
deg$(\tilde{T})=l$ is mapped to $T' = \tilde{T'}\cdot (m x'\partial' -
\alpha_1)(m x'\partial' - \alpha_2) \cdots (m x'\partial' -
\alpha_k)$ with deg$(\tilde{T'})=l/m$. This leads directly to the
proof of statement (2) of the Proposition. \hfill $\Box$
\end{enumerate}
The classification of second order operators $T_2 \in \D_V$ is similar
as before. A difference is that $|I^{(m)}| = |I^{(-m)}|$, so that
there is always a sort of symmetry in $T_2$. Here we give one
example, the generic 2-chain case. Suppose $|I^{(m)}| = 2$, so $I$ has
the structure $I= \{ i_1,i_3,\ldots,i_ {2r-1}\} \cup \{
i_2,i_4,\ldots,i_ {2s}\}$ with
\[ i_3=i_1-m,\ i_5=i_3-m,\ \ldots ,\ i_{2r-1}=i_{2r-3}-m \]
and
\[ i_4=i_2-m,\ i_6=i_4-m,\ \ldots ,\ i_{2s}=i_{2s-2}-m \]
If for all $l>0, l \neq m$ there holds $|I^{(l)}| > 2$, then the most
general $2^{nd}$ order operator $T_2 \in \D_V$ takes the form:
\[
T_2 =  \alpha_1x^m(x\partial-i_1)(x\partial-i_2)
     + \alpha_2x^2\partial^2 + \alpha_3x\partial + \alpha_4
     + \alpha_5x^{-m}(x\partial-i_{2r-1})(x\partial-i_{2s})
\]
Note the symmetry in $T_2$: there are as many terms of degree $m$ as
terms of degree $-m$.

{\bf Final remark.} In the case that $V = \langle 1,x,x^2,\ldots,x^{n}
\rangle$, the algebra $\D_V$ is essentially generated by
$\{\partial,2x\partial-n,x^2\partial-nx\}$, see \cite{T2}. These 3
elements form the Lie algebra $sl_2(\C)$, and hence $\D_V$ is
essentially some representation of the universal enveloping algebra of
$sl_2(\C)$. We were not able to find a similar structure in $\D_V$ for
general $V$, even for the simplest case $V=\langle
1,x,x^3 \rangle$. Particularly,
one can show that for the space  $V=\langle
1,x,x^3 \rangle$ the algebra $\D_V$ is an infinite-dimensional,
finite-generated algebra. It is defined by 11 generators, which are
differential operators of first-, second- and third order, and
the commutator of any two of them is expressible as an {\it ordered}
 cubic polynomial in these generators.
\\[\baselineskip]
{\bf Acknowledgement}\\
The second author (A.T.) wishes to thank M.A. Shubin for useful
discussion.


\newpage
\begin{thebibliography}{9}
\bibitem{B} S. Bochner, {\it \"Uber Sturm-Liouvillesche
   Polynomsysteme}, Math.
   Z. {\bf 29}, p. 730-736 (1929).
\bibitem{L} L.L. Littlejohn, {\it Orthogonal polynomial solutions to
   ordinary differential equations}, in ``Orthogonal polynomials and their
   applications'', eds. M. Alfaro e.a., Lecture Notes in Mathematics
   1329, Springer-Verlag, Berlin (1970).
\bibitem{T1}
   A.V. Turbiner, {\it Quasi-exactly solvable problems and $sl(2)$
   algebra},
   Comm. Math. Phys. {\bf 118}, p. 467-474 (1988).
\bibitem{T2}
   A.V. Turbiner, {\it Lie-algebras and polynomials},
   Journ. Phys. {\bf A25}
   L1087-L1093 (1992); \\ {\it Lie-algebraic approach to the
   theory of polynomial solutions. I. Ordinary differential
   equations and finite-difference equations in one variable},
   preprint CPT-91/P.2679 (1991);
   \\ {\it Lie algebras and linear operators with invariant subspace},
   preprint I.H.E.S./P/92/95 (December 1992 and June 1993
   (corrected version)), to appear in Lie Algebras, Cohomologies
   and New Findings in Quantum Mechanics, AMS Contemporary Math,
   N. Kamran and P. Olver (eds.), AMS, 1993
\end{thebibliography}
\end{document}